\documentclass[]{amsart}
\usepackage[textwidth=14cm,hcentering]{geometry}
\usepackage{amssymb}
\usepackage{amsthm}
\usepackage{amsmath,amscd}
\usepackage[mathscr]{euscript}
\usepackage[all]{xy}
\usepackage{float}
\usepackage[utf8]{inputenc}
\usepackage{lmodern}
\usepackage[T1]{fontenc}
\usepackage[colorlinks=true,linkcolor=red,citecolor=blue]{hyperref}
\usepackage{cleveref}
\usepackage{mathtools}
\usepackage{tikz}
\usepackage{bbm}
\usepackage{adjustbox}
\usepackage{todonotes}
\usepackage{pgf}
\usepackage{hieroglf}

\usetikzlibrary{cd}
\usetikzlibrary{arrows,positioning}

\usepackage{amssymb}
\usepackage{amsthm}
\usepackage{amsmath,amscd}
\usepackage[mathscr]{euscript}
\usepackage[all]{xy}
\usepackage{float}
\usepackage[utf8]{inputenc}
\usepackage{lmodern}
\usepackage[T1]{fontenc}
\usepackage[textwidth=14cm,hcentering]{geometry}
\usepackage[colorlinks=true,linkcolor=red,citecolor=blue]{hyperref}
\usepackage{mathtools}
\usepackage{tikz}
\usepackage{bbm}
\usepackage{adjustbox}
\usepackage{todonotes}
\usepackage{pgf}
\usetikzlibrary{cd}
\usetikzlibrary{arrows,positioning}

\usepackage{hyperref}
\usepackage{cleveref}
\theoremstyle{definition}
\newtheorem{thm}{Theorem}[section]

\newtheorem{rmq}[thm]{Remark}
\newtheorem{prop}[thm]{Proposition}
\newtheorem{defn}[thm]{Definition}
\newtheorem{lemme}[thm]{Lemma}
\newtheorem{coroll}[thm]{Corollary}

\newtheorem{thmalpha}{Theorem} 

\newcommand{\xto}[1]{\xrightarrow{#1}}

\newcommand{\bb}[1]{\mathbb{#1}}
\newcommand{\mc}[1]{\mathcal{#1}}

\newcommand{\mbf}[1]{\mathbf{#1}}
\newcommand{\mrm}[1]{\mathrm{#1}}

\newcommand{\bbZ}{\mathbb{Z}}

\newcommand{\ol}{\overline}

\usepackage{mathtools}
\DeclareUnicodeCharacter{202F}{ }

\usepackage{diagbox}

\usepackage{lineno}

\title{Tame homotopy theory and the Goodwillie tower of the identity}
\author{Nicolas Guès}

\crefname{thmalpha}{Theorem}{Theorems}
\newtheorem{coralpha}[thmalpha]{Corollary}
\crefname{coralpha}{Corollary}{Corollaries}
\begin{document}
	\begin{abstract}
	The purpose of this note is to provide a modern treatment of Dwyer's tame homotopy theory and extend it to non-simply connected spaces, exploiting a connection with the work of Heuts on the Goodwillie tower of the identity. As an application, we construct spectral Lie algebraic models for Hess's \emph{$r$-mild} homotopy theory. 
	\end{abstract}
	\maketitle
	\tableofcontents

	\section{Introduction}\label{section:intro}

	In 1979, Dwyer \cite{dwyer79} extended Quillen's approach to rational homotopy theory to a finer notion of weak equivalence. Fixing a connectivity $r \geq 3$, we say an $r$-tame weak equivalence is a continuous map $f: X \to Y$ between $(r-1)$-connected spaces such that for each $n \geq r$, $\pi_n f$ is an isomorphism after inverting each prime $p$ such that $r+2p-3\leq n$. The localization of $(r-1)$-connected pointed spaces with respect to $r$-tame weak equivalences is a reflective subcategory $\mc S_*^{r\text{-tame}}$ consisting of spaces $X$ such that multiplication by $p$ in $\pi_n X$ is invertible for each $p$ in the tame range $r+2p-3 \leq n$. Similarly, define an $r$-\emph{tame} dg Lie algebra to be an {$(r-1)$-reduced dg Lie algebra $L$ over $\bb Z$, i.e. $L_i = 0$ for $i < r-1$, such that $H_n L$ is a $\bb Z[1/p]$-module for each $p$ such that $r+2p-3 \leq n+1$ (in Quillen's convention $H_n L$ corresponds to $\pi_{n+1}$ of the space)}. Reframed in modern language, Dwyer proves the following equivalence of $\infty$-categories.
	\begin{thm}\cite{dwyer79}
		For $r\geq3$, there is an equivalence
		$$\mc S_*^{r\text{-tame}} \simeq \mrm{DGL}_{\bb Z}^{r\text{-tame}}$$
	\end{thm}
	where $\mrm{DGL}_{\bb Z}^{r\text{-tame}}$ is the category of DG Lie algebras over $\bb Z$ localized at $r$-tame equivalences.
	Similarly, $r$-tame homotopy types can be described by an algebraic homotopy category of tame commutative dg-algebras \cite{scheererReport88, berlinpaper, cenkl_porterCubical}, which we will not explore in this note.
 One of the objectives of this paper is to give a modern treatment of tame homotopy theory, and along the way, extend the range of its application. The first result is the following. We say a $2$-tame equivalence between simply connected spaces is a map $f$ inducing an isomorphism on $\pi_k$ after inverting $(k-1)!$. 
 
 \begin{thmalpha}\label[thmalpha]{2tame}
 	There is an equivalence of $\infty$-categories between simply connected spaces modulo $2$-tame equivalences and $1$-connected algebras over the integral Lie operad $\bb L(\bb Z)=H\bb Z\otimes\bb L$ whose $\pi_k$ is a $\bb Z[\frac1{(k-1)!}]$-module for every $k$:
	$$\mc S_*^{2\text{-tame}} \simeq \mrm{Alg}(\bb L(\bb Z))^{2\text{-tame}}$$
\end{thmalpha}
 This recovers a ``Koszul dual'' version of Scheerer's commutative algebra models constructed in \cite{scheererReport88}, see his Remark $2$. That the theory could be extended to simply connected spaces was already suggested in Dwyer's original paper \cite[Section 1.5]{dwyer79}: ``It seems clear that these problems could be avoided by restricting attention to ring systems with more divisibility than tame ring systems, for instance, to systems $R_*$ such that $\frac 1 p \in R_{p-2} (p\geq2)$''. Another generalization is Hess's \emph{mild} homotopy theory. For $r \geq 2$ a parameter, define an $r$-\emph{mild} homotopy equivalence \cite{hess93} as a map $f$ between $(r-1)$-connected spaces inducing an isomorphism on $\pi_k$ up to $k=p(r-1)$ after localizing at $p$. Note that this notion is finer than a tame equivalence (if a map is an $r$-mild equivalence then it is an $r$-tame equivalence). The second theorem of this note is that we can extend the scope of Dwyer's theory to \emph{mild} homotopy theory, at the cost of passing from $\bb Z$ to $\bb S$. 
\begin{thmalpha}\label[thmalpha]{mildmodel}
	 	For any $r \geq 2$, the category of $(r-1)$-connected spaces modulo $r$-mild equivalences is equivalent to the category of $(r-1)$-connected spectral Lie algebras modulo $r$-mild homotopy equivalences:
	$$\mc S_*^{r\text{-mild}} \simeq \mrm{Alg}(\bb L)^{r\text{-mild}}$$
\end{thmalpha}

The main idea of this paper can be summed up in the following slogan: in Dwyer's tame homotopy theory, one recovers from a space $X$ the collection of its Postnikov truncations $t_{n} X$ suitably localized (inverting primes up to the tame bound), together with coherences. A more fundamental and, we believe, more powerful approach is to replace the $n$-th Postnikov truncations with the $n$-th \emph{Goodwillie approximation} $\mc P_n X$ of the space $X$, suitably localized.  One can hence define an abstract $\infty$-category remembering only the collection of these localized Goodwillie approximations (which is not, contrary to tame homotopy theory, a localization of spaces). We call this $\bb S$-tame homotopy theory, and construct an equivalence with an $\infty$-category of algebraic objects that we call \emph{spectral Lazard algebras}. {\Cref{2tame,mildmodel} are obtained by truncating the following structural result.}
{
\begin{thmalpha}\label[thmalpha]{Stamemodel}
	There is an equivalence of $\infty$-categories
	$$\mc{FT} \simeq \mrm{LazAlg}(\bb S)$$
	between the $\infty$-category of $\bb S$-tame spaces (\Cref{formaltamespace}) and the $\infty$-category of spectral Lazard algebras (\Cref{def:spectrallazard}).
\end{thmalpha}
}
As the Goodwillie tower of the identity converges for simply connected spaces, the $\bb S$-tame homotopy type of a space allows one to recover its tame homotopy type, in Dwyer's sense, and thus the theorems above. An advantage of this notion, contrary to Dwyer's, is that it makes sense also for non-simply connected spaces.  This approach allows us to suggest a generalization of tame homotopy theory to non-simply connected spaces, embodied in the following {corollary of \Cref{Stamemodel}}, made more precise in \cref{section:pi1}.
\begin{coralpha}\label[coralpha]{pi1lazard}
	There is a map 
	$$\mc L:  \mc S_*^{\geq 1} \to \mrm {LazAlg}(\bb S)$$ 
	from pointed connected spaces to spectral Lazard algebras, that remembers:
	\begin{enumerate}
		\item {the Lie models of the $\mc P_n X[\frac 1 {n!}]$ given by \Cref{modelgoodwillie}, which are the stages of $\mc L X$, for all $n \in \bb N$;}
		\item the Lazard completion (see \cite{lazard54}) of the fundamental group of $X$. 
	\end{enumerate}
\end{coralpha}

\subsection{Outline}

In \Cref{section:lazard,section:Stamespaces} we deal with $\bb S$-tame homotopy theory and show how it is related to an algebraic object that we call a \emph{spectral Lazard algebra}, whose name is inspired by the classical notion of a Lazard algebra \cite{lazard54}. In particular, the {$\pi_1$} of a spectral Lazard algebra is a Lazard algebra in Lazard's original sense. In \Cref{section:mildandtame} we show how after truncating, this theory gives easily \Cref{2tame} and \Cref{mildmodel}. Finally in \Cref{section:pi1} we explore what $\bb S$-tame homotopy theory says on non-simply connected spaces.  All of these constructions will be further applied to construct algebraic models of configuration spaces and spaces of knots in \cite{gues26}.

\subsection{Conventions}
We will say a map $X \to Y$ between spaces, spectra or chain complexes is $n$-connected if its homotopy fiber is $(n-1)$-connected. {A spectrum $E$ is $k$-connective if $\pi_i E = 0$ for $i < k$. We write $t_n$ for the Postnikov truncation killing the homotopy groups in degrees $> n$, and $\tau_n \mc O$ for the truncation of an operad $\mc O$ killing the operations of arity $> n$. We write $\bb L$ for the spectral Lie operad, i.e. Ching's operad $\partial_*(\mrm{id})$ of derivatives of the identity, and $\bb L(\bb Z) = H\bb Z \otimes \bb L$. With this convention there is no shift: the Lie model $\mc L(X)$ of a space $X$ satisfies $\pi_k \mc L(X) \cong \pi_k X$ (\Cref{modelgoodwillie}), whereas in Quillen's and Dwyer's dg Lie models the homology in degree $k$ corresponds to $\pi_{k+1}$.}
\subsection{Acknowledgments}
I thank my advisor, Geoffroy Horel, for many insightful discussions and indispensable guidance. I thank Gregory Arone for help in the proof of \Cref{pi1goodwillie}.

\subsection{Use of generative AI}
The author acknowledges the use of large language models for help in the writing and proofreading of this manuscript. No generative AI tool was used in the mathematical creative process.
	\section{Lazard algebras and spectral Lazard algebras}\label{section:lazard}
	
	A category of groups, which we will call Lazard groups following \cite{cenklporterlazard}, was introduced by Lazard \cite{lazard54}. It is the natural family of groups such that the Baker--Campbell--Hausdorff formula makes sense. 
	
	To each group $G$ we can associate its Lazard completion $\ol G$ which is the initial Lazard group to which $G$ maps. The Lazard completion is a stronger invariant than the Malcev completion as it contains partial information about torsion. We refer to \cite[Section 1]{cenklporterlazard}, \cite[Chap.~II, \S 4]{lazard54} for a detailed treatment of this topic. We define next what we need to deal with Lazard groups and Lazard algebras.
	For $L$ a Lie algebra over $\bbZ$, a descending filtration $F$ on $L$ is a decreasing sequence $$L = F_1 L \supset F_2 L \supset \cdots$$
	such that $[F_i L, F_j L] \subset F_{i+j}L$. Thus the associated graded $\mathrm{gr}(L) = \bigoplus F_i L/F_{i+1} L$ has a Lie algebra structure. A descending filtration of a group $G$ is  an analogous decreasing sequence $F_\bullet G$ with  $[F_i G, F_j G] \subset F_{i+j} G$.
	For $i$ an integer, write $Q_i = \bbZ[\frac{1}{i!}]$.
	
	\begin{defn}
		A \emph{Lazard algebra} is a Lie algebra $L$ over $\bb Z$ endowed with a filtration $F_\bullet L$ satisfying:  
		\begin{enumerate}
			\item For each $i \geq1$, $F_i L$ is a $Q_i$-module.
			\item For each $i \geq 1$, the nilpotent Lie algebra $\mathrm{gr}(L/F_{i+1} L) \otimes Q_i$ is generated by $(L/F_2 L)\otimes Q_i$. 
			\item $L$ is complete with respect to the filtration $F_\bullet L$.  
		\end{enumerate}
	\end{defn}
	\begin{defn}
		A \emph{Lazard group} is a group $G$ endowed with a descending filtration $F_\bullet G$ satisfying:  
		\begin{enumerate}
			\item For each $i \geq1$, $F_i G/F_{i+1} G$ is a $Q_i$-module.
			\item For each $i \geq 1$, the nilpotent Lie algebra $\mathrm{gr}(G/F_{i+1} G) \otimes Q_i$ is generated by $(G/F_2 G)\otimes Q_i$. 
			\item $G$ is complete with respect to the filtration $F_\bullet G$.  
		\end{enumerate}
	\end{defn}
	{Lazard algebras (resp. Lazard groups) form a category, whose morphisms are the homomorphisms of Lie algebras (resp. of groups) preserving the filtrations.} 
	Lazard \cite{lazard54} showed that the BCH formula induces an equivalence between the category of Lazard algebras and the category of Lazard groups.
	\begin{defn}\cite[\S 1]{cenklporterlazard}, see also \cite[Chap.~II, Thm.~3.10]{lazard54}
		For each group $G$ there exists a Lazard group $\ol G$  together with a map $G \xto g \ol G$ such that if $G$ is filtered by its lower central series subgroups then $g$ is filtration preserving and the induced map 
		$$\bigoplus (F_iG/F_{i+1} G) \otimes Q_i \to \mathrm{gr}(\ol G)$$ is an isomorphism. $\ol G$ is called the Lazard completion of the group $G$.
	\end{defn}
	There exists an analogous Lazard completion of a Lie algebra over $\bb Z$ defined by a universal property.
	\begin{rmq}
		Note that if the rings $Q_i$ are replaced by the rationals $\bb Q$ in these definitions we recover the notions of Malcev complete Lie algebras and Malcev complete groups, as well as the Malcev completion.
	\end{rmq}
	
	\begin{prop}\label[prop]{equivalencelazarddefs}
		The category of Lazard algebras is equivalent to the limit (in the $(\infty,1)$-category of $(\infty,1)$-categories) of the following diagram of $1$-categories: 
		\[\begin{tikzcd}[row sep=large, column sep=-30]
			{\mathrm{Alg}(\tau_{1} \mathrm{Lie}(\mathbb Z))} &&
			{\mathrm{Alg}(\tau_{2} \mathrm{Lie}(\mathbb Z))\bigl[\tfrac12\bigr]} &&
			{\mathrm{Alg}(\tau_{3} \mathrm{Lie}(\mathbb Z))\bigl[\tfrac16\bigr]} &&
			\cdots \\
			&
			{\mathrm{Alg}(\tau_{1} \mathrm{Lie}(\mathbb Z))\bigl[\tfrac12\bigr]} &&
			{\mathrm{Alg}(\tau_{2} \mathrm{Lie}(\mathbb Z))\bigl[\tfrac16\bigr]} &&
			{\mathrm{Alg}(\tau_{3} \mathrm{Lie}(\mathbb Z))\bigl[\tfrac1{24}\bigr]}
			\arrow[from=1-1, to=2-2]
			\arrow[from=1-3, to=2-2]
			\arrow[from=1-3, to=2-4]
			\arrow[from=1-5, to=2-4]
			\arrow[from=1-5, to=2-6]
			\arrow[from=1-7, to=2-6]
		\end{tikzcd}\]
	\end{prop}
	Here $\mrm{Lie}(\bb Z)$ denotes the ($1$-categorical) Lie operad and $\tau_n$ its truncation killing operations of arity $> n$. Hence the category of algebras over $\tau_n \mrm{Lie}(\bb Z)$ is simply the category of $n$-nilpotent Lie algebras. 
\begin{proof}
	For now call $M$ the limit category above and $\mrm{Laz}$ the category of Lazard algebras. We will construct functors $M \leftrightarrows \mrm{Laz}$. An object of $M$ is a sequence $(L_n)_{n\geq1}$, where $L_n$ is an $n$-nilpotent Lie algebra over $Q_n$, together with isomorphisms
	\[\theta_n\colon L_{n+1}/\Gamma^{n+1}L_{n+1}\xto{\simeq}L_n\otimes Q_{n+1}.\]
	Iterating, they induce isomorphisms
	\begin{equation}\label{eq:thetaitere}
		L_k/\Gamma^jL_k\cong L_{j-1}\otimes Q_k\qquad(j\leq k).
	\end{equation}
	We will use the following observation: if $0\to K\to A\to B\to0$ is an exact sequence of
	abelian groups and $K$ is a $Q_n$-module, then
	\begin{equation}\label{eq:carrecartesien}
		A\cong(A\otimes Q_n)\times_{B\otimes Q_n}B,
	\end{equation}
	since $A\to B$ and $A\otimes Q_n\to B\otimes Q_n$ are surjective with the same kernel $K$.
	
	\textbf{From $\mrm{Laz}$ to $M$.} Let $L$ be a Lazard algebra, and set $E_n=L/F_{n+1}L$ and
	$L_n=E_n\otimes Q_n$. Since $\Gamma^{n+1}L\subset F_{n+1}L$, this is an $n$-nilpotent Lie
	algebra over $Q_n$. By condition (2) for $i=n+1$, the associated graded of
	$E_{n+1}\otimes Q_{n+1}$ is generated in degree $1$; by descending induction on $j$, this
	implies that $F_j(E_{n+1}\otimes Q_{n+1})=\Gamma^j(E_{n+1}\otimes Q_{n+1})$ for all $j$. For
	$j=n+1$ we get
	\[L_{n+1}/\Gamma^{n+1}L_{n+1}=(E_{n+1}/F_{n+1}E_{n+1})\otimes Q_{n+1}=E_n\otimes Q_{n+1}
	=L_n\otimes Q_{n+1},\]
	which defines $\theta_n$. This construction is functorial.
	
	\textbf{From $M$ to $\mrm{Laz}$.} Let $(L_n,\theta_n)$ be an object of $M$. Define Lie rings
	$E_n$ by $E_1=L_1$ and
	\[E_{n+1}=L_{n+1}\times_{L_n\otimes Q_{n+1}}E_n,\]
	where $E_n\to L_n$ is the first projection. Explicitly, $E_n$ is the set of sequences
	$(x_n,\dots,x_1)$ with $x_k\in L_k$, such that $x_{k+1}$ and $x_k$ have the same image in
	$L_k\otimes Q_{k+1}$. The projection $E_{n+1}\to E_n$ is surjective with kernel
	$\Gamma^{n+1}L_{n+1}$, and $E_n\otimes Q_n\cong L_n$ by induction, as localization commutes
	with fiber products. Let $L=\lim_nE_n$ and $F_{n+1}L=\ker(L\to E_n)$. Then $L$ is complete,
	and $F_iL/F_{i+1}L=\Gamma^iL_i$ is a $Q_i$-module, so that $F_iL$ is a $Q_i$-module. By
	\eqref{eq:thetaitere}, an element $x$ lies in $F_jL$ if and only if $x_k=0$ for $k<j$ and
	$x_k\in\Gamma^jL_k$ for $k\geq j$. Hence $[F_iL,F_jL]\subset F_{i+j}L$, and the filtration
	induced on $(L/F_{i+1}L)\otimes Q_i\cong L_i$ is the lower central series, whose associated
	graded is generated in degree $1$. So $L$ is a Lazard algebra, functorially in
	$(L_n,\theta_n)$.
	
Now let us show that the two functors are inverse to each other. Starting from $(L_n,\theta_n)$, and calling $L$ the Lazard algebra constructed before, we have
	$(L/F_{n+1}L)\otimes Q_n=E_n\otimes Q_n\cong L_n$, compatibly with the $\theta_n$. Conversely,
	let $L$ be a Lazard algebra and $E_n=L/F_{n+1}L$, $L_n=E_n\otimes Q_n$ as above. Applying
	\eqref{eq:carrecartesien} to the extension $0\to F_{n+1}L/F_{n+2}L\to E_{n+1}\to E_n\to0$,
	whose kernel is a $Q_{n+1}$-module, we get
	\[E_{n+1}\cong(E_{n+1}\otimes Q_{n+1})\times_{E_n\otimes Q_{n+1}}E_n
	=L_{n+1}\times_{L_n\otimes Q_{n+1}}E_n.\]
	Thus the Lie rings $E_n$ are the ones built from $(L_n,\theta_n)$, and
	$L\cong\lim_nE_n$ with its filtration, since $L$ is complete.
\end{proof}
	This result motivates the following definition. 
	\begin{defn}\label[defn]{def:spectrallazard}
		The category {$\mrm{LazAlg}(\bb S)$} of spectral Lazard algebras is the following limit of $\infty$-categories. 
		\[\begin{tikzcd}[row sep=normal, column sep=-30]
			{\mathrm{Alg}(\tau_{1} \bb L)^{\geq 1}} &&
			{\mathrm{Alg}(\tau_{2} \bb L)^{\geq 1}\bigl[\tfrac12\bigr]} &&
			{\mathrm{Alg}(\tau_{3} \bb L)^{\geq 1}\bigl[\tfrac16\bigr]} &&
			\cdots \\
			&
			{\mathrm{Alg}(\tau_{1} \bb L)^{\geq 1}\bigl[\tfrac12\bigr]} &&
			{\mathrm{Alg}(\tau_{2} \bb L)^{\geq 1}\bigl[\tfrac16\bigr]} &&
			{\mathrm{Alg}(\tau_{3} \bb L)^{\geq 1}\bigl[\tfrac1{24}\bigr]}
			\arrow[from=1-1, to=2-2]
			\arrow[from=1-3, to=2-2]
			\arrow[from=1-3, to=2-4]
			\arrow[from=1-5, to=2-4]
			\arrow[from=1-5, to=2-6]
			\arrow[from=1-7, to=2-6]
		\end{tikzcd}\]
		
	\end{defn}
	{Here $\mrm{Alg}(\tau_k \bb L)$ denotes the $\infty$-category of algebras in spectra over $\tau_k\bb L$, whereas in \Cref{equivalencelazarddefs} $\mrm{Alg}$ denotes ordinary categories of algebras. The superscript $\geq 1$ means that the underlying spectrum is $1$-connective. With our conventions (no shift), the bracket of an $\bb L$-algebra $A$ induces maps $\pi_a A \otimes \pi_b A \to \pi_{a+b-1} A$; in particular $\pi_1 A$ is a Lie algebra, in degree $1$ rather than $0$ as in \Cref{equivalencelazarddefs}.}
	\section{The $\infty$-category of $\bb S$-tame spaces}\label{section:Stamespaces}

The upshot of this section is that the category of $\bb S$-tame spaces can be modeled by the spectral Lazard algebras of \Cref{section:lazard}.
In \cite{heuts2021}, Heuts constructs, for any pointed presentable $\infty$-category $\mc C$, a tower of adjunctions (we write only the left adjoints)
$$\mc C \to \cdots \to \mc P_n \mc C \to \cdots \to \mc P_0 \mc C = *$$
where $\mc P_n \mc C$ is in some precise sense the ``universal $n$-excisive approximation of $\mc C$''. In particular, the $n$-th stage of the Goodwillie tower of the identity of $\mc C$ can be expressed as the unit of the adjunction: 
$$\mc C \xrightleftharpoons[\Omega^\infty_n]{\Sigma^\infty_n} \mc P_n \mc C$$
These $\mc P_n$ approximations are defined using the machinery of $\infty$-operads. In this section, this notion will be thought of as Lurie's definition in \cite[Chapter 2, Definition 2.1.1.10]{Lurie2015} to remain coherent with Heuts' work. 
Let us recall the notion of the $n$-truncation of an $\infty$-operad $\mc O$ in the sense of Heuts. Informally, $\tau_n \mc O$ is an $\infty$-operad whose underlying symmetric sequence is defined by
\[
\tau_n \mathcal{O}(k) =
\begin{cases}
	\mathcal{O}(k) & \text{if } 1 \leq k \leq n, \\
	* & \text{if } k > n.
\end{cases}
\]
The operad structure on $\tau_n \mc O$ is inherited in the obvious way from $\mc O$, simply setting operations of arity greater than $n$ to zero. There is an obvious truncation map $\mc O \to \tau_n \mc O$; see \cite[Section 4.1]{heuts2021} for the more rigorous construction. To construct $\bb S$-tame homotopy theory, we need a bit more background on $\infty$-operads. What follows is a short summary of some tools developed in \cite{heuts2021}. 
	
	Throughout, $\mc O^{\otimes}$ denotes a nonunital stable $\infty$-operad in the sense of \cite{heuts2021}, with underlying stable $\infty$-category $\mc O$. Such an operad comes with $k$-fold tensor products $\otimes^k \colon \mc O^k \to \mc O$ for $k \geq 1$; for $X \in \mc O$ we write $X^{\otimes k}$ for $\otimes^k(X, \dots, X)$, which carries an action of $\Sigma_k$. The example we care about is $\mrm{Sp}^{\otimes}$, for which $X^{\otimes k} = X^{\wedge k}$.

	\begin{defn}\label[defn]{def:coalgebra}
		A \emph{coalgebra} in $\mc O^{\otimes}$ is an object $X \in \mc O$ together with a coherent family of comultiplications
		\[
		\delta_k \colon X \longrightarrow (X^{\otimes k})^{h\Sigma_k}, \qquad k \geq 2,
		\]
		see \cite[Def.~4.14]{heuts2021} for the precise definition. We write $\mrm{coAlg}(\mc O^{\otimes})$ for the $\infty$-category of coalgebras and $\mrm{coAlg}^{\mrm{ind}}(\mc O^{\otimes})$ for its ind-completed variant \cite[Def.~4.21]{heuts2021}. Since $\mc O$ is stable, every $X \in \mc O$ carries the \emph{trivial} coalgebra structure, for which all the $\delta_k$ are null.
	\end{defn}

	\begin{rmq}\label[rmq]{rmq:stabilization}
		The functor $\Sigma^\infty_n \colon \mc C \to \mc P_n \mc C$ induces an equivalence on stabilizations \cite[Lem.~2.11]{heuts2021}. For $\mc C = \mc S_*$ this identifies $\mrm{Sp}(\mc P_n \mc S_*) \simeq \mrm{Sp}$, and under this identification the spectrum underlying $\Sigma^\infty_n X$ is $\Sigma^\infty X$.
	\end{rmq}

	Fix $n \geq 1$. The truncation $\tau_n \mc O^{\otimes}$ has the same underlying $\infty$-category $\mc O$ but different tensor products, which we denote $\odot^k$. They come with natural transformations $\otimes^k \to \odot^k$, which are equivalences for $k \leq n$. For $k > n$, $\odot^k$ is a limit of iterated tensor products, indexed by the nerve $\mrm N\mbf{Part}_n(k)$ of a poset of chains of partitions of $\{1, \dots, k\}$, which records the ways of building a $k$-ary operation out of operations of arity $\leq n$ \cite[\S 4.1, Prop.~4.10]{heuts2021}. For $\mc O^{\otimes} = \mrm{Sp}^{\otimes}$ all these iterated tensor products are copies of $X^{\wedge k}$. The only feature of this description we need is that the nondegenerate simplices of $\mrm N\mbf{Part}_{k-1}(k)$ have dimension at most $k-3$ (see the proof of \cite[Lemma 6.17]{heuts2021}).

	\begin{lemme}\label[lemme]{lemme:liftcriterion}
		Let $k \geq 2$ and let $X \in \mc O$ be endowed with the trivial $\tau_{k-1}\mc O^{\otimes}$-coalgebra structure. Extensions of this structure to a $\tau_k \mc O^{\otimes}$-coalgebra structure on $X$ are lifts of the zero map $X \to (X^{\odot k})^{h\Sigma_k}$ along $(X^{\otimes k})^{h\Sigma_k} \to (X^{\odot k})^{h\Sigma_k}$ \cite[Lemma~4.23]{heuts2021}. If the Tate construction of $\Sigma_k$ vanishes in $\mc O$, for instance if $k!$ is invertible, the norm map identifies homotopy fixed points with homotopy orbits, and these extensions are classified by the space $\mrm{Map}_{\mc O}(X, F)$, where $F = \mrm{fib}\bigl(X^{\otimes k} \to X^{\odot k}\bigr)_{h\Sigma_k}$ \cite[Cor.~4.32]{heuts2021}. The trivial extension corresponds to the zero map. In particular, if every map $X \to F$ is null, every such extension is equivalent to the trivial one.
	\end{lemme}
	
	\begin{thm}\cite[Corollaries 1.26 and 6.15]{heuts2021}.
		After inverting $(n-1)!$, there is an equivalence of $\infty$-categories
		$$\mc P_{n-1} \mc S_* \simeq \mrm{coAlg}^{\mrm{ind}}(\tau_{n-1}\mrm{Sp}^\otimes)$$
	\end{thm}

	Fix an integer $n$. Consider the $(n-1)$-th Goodwillie approximation of the category of pointed spaces $\mc P_{n-1} \mc S_*$. 
	\begin{defn}
		Call $\mc P_{n-1} \mc S_*^{\geq 1}$ the subcategory of $\mc P_{n-1} \mc S_*$ generated under colimits by $\Sigma_{n-1}^{\infty} S^1$. 
	\end{defn}
We recall that for any $n$-truncated operad $\mc O$ in spectra there exists a functor
$$\mbf{bar}_{\mc O}\colon \mrm{Alg}(\mc O) \longrightarrow \mrm{coAlg}_{\mrm{dp}}(B(\mc O))$$
to divided power coalgebras \cite[\S 6.1, Prop.~6.1]{heuts2021} over the bar construction $B(\mc O)$, which is a cooperad. On underlying spectra, $\mbf{bar}_{\mc O}(A)$ is the spectrum of derived indecomposables $\mathbf 1 \circ_{\mc O} A$, and $\mbf{bar}_{\mc O} \circ \mrm{free}_{\mc O}$ is the trivial coalgebra functor. This is a spectral, operadic analogue of Quillen's correspondence between dg Lie algebras and cocommutative dg coalgebras. With this analogy in mind, $\mbf{bar}$ must be thought of as the Chevalley--Eilenberg complex. The functor $\mbf{bar}_{\mc O}$ is fully faithful, with essential image the conilpotent coalgebras \cite[Prop.~6.9]{heuts2021}. For $\mc O = \tau_{n-1}\bb L$, the bar construction $B(\tau_{n-1}\bb L) \simeq \tau_{n-1}B(\bb L)$ is the $(n-1)$-truncated cocommutative cooperad, whose coalgebras are the coalgebras in $\tau_{n-1}\mrm{Sp}^{\otimes}$, and divided power coalgebras are ordinary coalgebras if the Tate constructions of the $\Sigma_k$ vanish for  $k \leq n-1$. Hence, after inverting $(n-1)!$, we obtain a fully faithful functor
$$\mbf{bar}_{\tau_{n-1}\bb L}\colon \mrm{Alg}(\tau_{n-1}\bb L) \longrightarrow \mrm{coAlg}^{\mrm{ind}}(\tau_{n-1}\mrm{Sp}^{\otimes}).$$
See \cite[\S 6.1 and Cor.~6.16]{heuts2021} for more details.
	\begin{lemme}
		The map $ \mc P_{n-1} \mc S_*^{\geq 1} \to \mrm{coAlg}^{\mrm {ind}}(\tau_{n-1} \mrm{Sp}^{\otimes})$ lands in the essential image of the functor $\mbf{bar}_{\tau_{n-1} \bb L}$.
	\end{lemme}
	\begin{proof}
		It is sufficient to prove it for $\Sigma^\infty_{n-1} S^1$ because the functor above
		preserves colimits. We adapt the proof of \cite[Lemma 6.17]{heuts2021}. We prove that
		$\Sigma^\infty_{n-1} S^1$ is a trivial coalgebra in $\tau_{n-1} \mrm{Sp}^{\otimes}$,
		by induction on $1 \leq k \leq n-2$.
		Here $\odot^{k+1}$ is now the $(k+1)$-fold tensor product defined by
		$\tau_k\mathrm{Sp}^{\otimes}$. We claim that any lift of $\Sigma^\infty S^1$ to
		$\tau_{k+1}\mathrm{Sp}^{\otimes}$ is uniquely determined and is the trivial one.
		Such a lift is a factorization of $\delta_{<k+1}$, which vanishes by the induction
		hypothesis, through
		$\bigl((\Sigma^\infty S^1)^{\otimes k+1}\bigr)^{h\Sigma_{k+1}} \to
		\bigl((\Sigma^\infty S^1)^{\odot k+1}\bigr)^{h\Sigma_{k+1}}$
		{(\Cref{lemme:liftcriterion})},
		so it suffices to prove that
		{the fiber $F$ of this map satisfies $\pi_1 F = 0$}.
		Since $k+1 \leq n-1$, the order of $\Sigma_{k+1}$ divides $(n-1)!$ and is invertible, so the norm map
		$(-)_{h\Sigma_{k+1}} \to (-)^{h\Sigma_{k+1}}$ is an equivalence and homotopy fixed
		points preserve connectivity. As $(\Sigma^\infty S^1)^{\otimes k+1}$ is
		$2$-connective for $k \geq 1$, so is
		$\bigl((\Sigma^\infty S^1)^{\otimes k+1}\bigr)^{h\Sigma_{k+1}}${. Moreover,
		$(\Sigma^\infty S^1)^{\odot k+1}$ is a limit of copies of the $(k+1)$-connective spectrum
		$(\Sigma^\infty S^1)^{\wedge k+1}$ indexed by $\mrm N\mbf{Part}_k(k+1)$, whose
		nondegenerate simplices have dimension at most $k-2$; it is therefore $3$-connective,
		and so is $\bigl((\Sigma^\infty S^1)^{\odot k+1}\bigr)^{h\Sigma_{k+1}}$. Hence $F$ is
		$2$-connective}, and the only lift of
		the zero map is the trivial one.
	\end{proof}
	\begin{prop}\label[prop]{modelgoodwillie}
		After inverting $(n-1)!$, there is an equivalence of $\infty$-categories
		$$\mc L: \mc P_{n-1} \mc S_*^{\geq 1} \simeq \mrm{Alg}(\tau_{n-1} \bb L)^{\geq 1}$$
		Moreover, for $X \in \mc P_{n-1} \mc S_*^{\geq 1}$, if we write $\pi_k X = \pi_0 \mrm{Map}_{\mc P_{n-1} \mc S_*}(\Sigma^{\infty}_{n-1}S^k, X)$, we have $\pi_k X \simeq \pi_k(\mc L(X))$. For $k=1$, this isomorphism is only set-theoretic while if $k \geq2$ it is an abelian group isomorphism.
	\end{prop}
	\begin{proof}
		By construction, the functor $\mc L$ admits a right adjoint $\mc R$. Let us prove that $ \mc L$ is fully faithful. We have an equivalence $\mc P_{n-1} \mc S_* \xto{\sim} \mrm{coAlg}^{\mrm{ind}}(\tau_{n-1} \mrm{Sp}^{\otimes})$. Restricting to $\mc P_{n-1} \mc S_*^{\geq 1}$ this map factors through $\mrm{Alg}(\tau_{n-1} \bb L)$ which we can identify as a full subcategory of $\mrm{coAlg}^{\mrm{ind}}(\tau_{n-1} \mrm{Sp}^{\otimes})$ by the bar construction and \cite[Cor.~6.16]{heuts2021}. Hence $\mc L$ is fully faithful. Now let us prove that $\mc L$ is essentially surjective. But $\mc L$ preserves colimits and its essential image contains connected free truncated Lie algebras (they are the images of the spheres $\Sigma^{\infty}_{n-1} S^k$ for $k \geq 1$). Any connected algebra $L$ in $\mrm{Alg}(\tau_{n-1} \bb L)$ is a colimit of such free algebras, for example using the bar resolution {and the fact that the free $\tau_{n-1}\bb L$-algebra on a $1$-connective spectrum is $1$-connective}. {Since $\mc L$ is fully faithful, such a diagram of free algebras lifts to a diagram in $\mc P_{n-1}\mc S_*^{\geq 1}$, and since $\mc L$ preserves colimits it sends the colimit of the lifted diagram to the colimit of the original one.} 
		{Let us now explain the identification of homotopy groups. Fix $k \geq 1$. As proved above, $\mc L(\Sigma^\infty_{n-1} S^k)$ is the free $\tau_{n-1}\bb L$-algebra $\mrm{free}(\Sigma^\infty S^k)$. Since $\mc L$ is fully faithful, we obtain natural bijections
		$$\pi_k X = \pi_0 \mrm{Map}(\Sigma^\infty_{n-1} S^k, X) \cong \pi_0 \mrm{Map}_{\mrm{Alg}}(\mrm{free}(\Sigma^\infty S^k), \mc L X) \cong \pi_0 \mrm{Map}_{\mrm{Sp}}(\Sigma^\infty S^k, U\mc L X) = \pi_k \mc L(X),$$
		where $U$ denotes the underlying spectrum. The group structure on $\pi_k X$ is induced by the pinch map $S^k \to S^k \vee S^k$. Since $\mc L$ preserves coproducts, the pinch map corresponds under $\mc L$ to a map of algebras
		$$\mrm{free}(\Sigma^\infty S^k) \longrightarrow \mrm{free}(\Sigma^\infty S^k \vee \Sigma^\infty S^k),$$
		which by adjunction is determined by a class in $\pi_k U\mrm{free}(\Sigma^\infty S^k \vee \Sigma^\infty S^k)$. Denote by $e_1, e_2$ the classes in this group given by the inclusions of the two copies of $\Sigma^\infty S^k$, and by $e$ the class of the map above. The arity $j$ summand of $\mrm{free}(e_1,e_2)$ is $(j(k-1)+1)$-connective. Hence for $k \geq 2$ only the arity $1$ summand contributes to $\pi_k\,\mrm{free}(e_1,e_2)$, which is therefore generated by $e_1$ and $e_2$, and compatibility with the two collapse maps $S^k \vee S^k \to S^k$ forces $e = e_1 + e_2$: This proves that the group law on $\pi_k X$ is given by addition in $\pi_k \mc L X$. }
	\end{proof}
	Note that  when $k=1$, the group isomorphism breaks down, because all arities contribute to $\pi_1\,\mrm{free}(e_1,e_2)$, which is the free $(n-1)$-nilpotent Lie algebra on $e_1$ and $e_2$, and the image of $e$ is a universal Lie series $e_1 + e_2 + (\text{brackets})$. The group law on $\pi_1 X$ thus in this case should correspond to a truncated Baker--Campbell--Hausdorff group law on $\pi_1 \mc L(X)$, and not to its addition.
	
	Truncating below $2p-3$ we get a model in Lie algebras over $\bb Z$.
	\begin{coroll}\label[coroll]{modelgoodwillieZ}
		After inverting $(p-1)!$ there is an equivalence of $\infty$-categories
		$$\mc P_{p-1} \mc S_*^{[1, 2p-3]} \simeq \mrm{Alg}(\tau_{p-1} \bb L(\bb Z))^{[1,2p-3]}$$
	\end{coroll}
	{Here $\mc P_{p-1}\mc S_*^{[1,2p-3]}$ denotes the full subcategory of $\mc P_{p-1}\mc S_*^{\geq 1}$ spanned by the objects $X$ with $\pi_k X = 0$ for $k > 2p-3$, and $\mrm{Alg}(\tau_{p-1}\bb L(\bb Z))^{[1,2p-3]}$ the $\infty$-category of algebras in spectra over the operad $\tau_{p-1}\bb L(\bb Z) = H\bb Z \otimes \tau_{p-1}\bb L$, equivalently of $\tau_{p-1}\bb L$-algebras in $H\bb Z$-modules, whose underlying spectrum has homotopy groups concentrated in degrees $[1,2p-3]$.}
	\begin{proof}
		{Throughout, $(p-1)!$ is inverted. By \Cref{modelgoodwillie} applied with $n = p$, which identifies homotopy groups, the left hand side is equivalent to $\mrm{Alg}(\tau_{p-1}\bb L)^{[1,2p-3]}$. We cannot truncate the operad $\tau_{p-1}\bb L$ directly: its terms are not connective, $\bb L(k)$ starting in degree $1-k$, and Postnikov truncations are only compatible with smash products on connective spectra. We therefore first shift. For an operad $\mc O$ in spectra, let $s\mc O$ be its operadic suspension, with $s\mc O(k) \simeq \Sigma^{k-1}\mc O(k)$, the symmetric group $\Sigma_k$ acting also on the suspension coordinates; the functor $A \mapsto \Sigma^{-1}A$ induces an equivalence $\mrm{Alg}(\mc O) \simeq \mrm{Alg}(s\mc O)$. Since $\bb L(k)$ is non-equivariantly a wedge of copies of $S^{1-k}$, the operad $s\tau_{p-1}\bb L$ is connective, its terms being non-equivariantly wedges of copies of $\bb S$, and
		$$\mrm{Alg}(\tau_{p-1}\bb L)^{[1,2p-3]} \simeq \mrm{Alg}(s\tau_{p-1}\bb L)^{[0,2p-4]}.$$
		On connective spectra, $t_{2p-4}$ is a symmetric monoidal localization \cite[\S 2.2.1]{Lurie2015}. Hence $t_{2p-4}s\tau_{p-1}\bb L$ is an operad in $(2p-4)$-truncated connective spectra, and since algebras with homotopy groups in $[0,2p-4]$ are local objects, $\mrm{Alg}(s\tau_{p-1}\bb L)^{[0,2p-4]} \simeq \mrm{Alg}(t_{2p-4}s\tau_{p-1}\bb L)^{[0,2p-4]}$. The same holds for the connective operad $H\bb Z \otimes s\tau_{p-1}\bb L$. For $0 < j \leq 2p-4$, the order of the finite group $\pi_j \bb S$ is only divisible by primes $q$ with $2q-3 \leq j$, hence $q < p$; therefore $t_{2p-4}\bb S \simeq H\bb Z$ once $(p-1)!$ is inverted. As the terms of $s\tau_{p-1}\bb L$ are non-equivariantly wedges of copies of $\bb S$, the unit $\bb S \to H\bb Z$ induces an equivalence of operads
		$$t_{2p-4}s\tau_{p-1}\bb L \xrightarrow{\ \sim\ } t_{2p-4}(H\bb Z \otimes s\tau_{p-1}\bb L).$$
		Putting this together and shifting back, we obtain
		$$\mrm{Alg}(\tau_{p-1}\bb L)^{[1,2p-3]} \simeq \mrm{Alg}(t_{2p-4}s\tau_{p-1}\bb L)^{[0,2p-4]} \simeq \mrm{Alg}(t_{2p-4}(H\bb Z\otimes s\tau_{p-1}\bb L))^{[0,2p-4]} \simeq \mrm{Alg}(\tau_{p-1}\bb L(\bb Z))^{[1,2p-3]}.$$}
	\end{proof}
This result can be put in perspective with \cite[Lemma 3.7.13]{wu22}.

	For $n\geq 1$, we will write $\mc P_n \mc S_*^{\geq 1}[\frac 1{k!}]$ to denote the localization of the category $\mc P_n \mc S_*^{\geq 1}$ at the equivalences defined by the local objects $X$ such that $\pi_* X$ is uniquely $k!$-divisible. 
	Informally, an $\bb S$-tame space is a sequence $(X_n)$ such that $X_n \in \mc P_n \mc S_*^{\geq 1}[\frac 1 {n!}]$ with a specified equivalence $X_n[\frac 1 {n+1}] \simeq \mc P_n X_{n+1}$.
	\begin{defn}\label[defn]{formaltamespace}
		The category $\mc {FT}$ of $\bb S$-tame spaces is the limit of the following diagram of $\infty$-categories: 
		\[\begin{tikzcd}[row sep=large, column sep=-10]
			{\mc P_1 \mc S_*^{\geq1}} && {\mc P_2 \mc S_*^{\geq1}[\frac 1 2]} && {\mc P_3 \mc S_*^{\geq1}[\frac 1 6]} && \cdots \\
			& {\mc P_1 \mc S_*^{\geq1}[\frac 1 2]} && {\mc P_2 \mc S_*^{\geq1}[\frac 1 6]} && {\mc P_3 \mc S_*^{\geq1}[\frac 1 {24}]}
			\arrow[from=1-1, to=2-2]
			\arrow[from=1-3, to=2-2]
			\arrow[from=1-3, to=2-4]
			\arrow[from=1-5, to=2-4]
			\arrow[from=1-5, to=2-6]
			\arrow[from=1-7, to=2-6]
		\end{tikzcd}\]
		
	\end{defn}
	\begin{rmq}
		One could think of another, maybe more natural diagram to take the limit of: 
		\[\begin{tikzcd}[row sep=small, column sep=-10]
			{\mc P_1 \mc S_*^{\geq1}} && {\mc P_2 \mc S_*^{\geq1}[\frac 1 2]} && {\mc P_3 \mc S_*^{\geq1}[\frac 1 6]} && \cdots \\
			& {\mc P_1 \mc S_*^{\geq1}[\frac 1 2]} && {\mc P_2 \mc S_*^{\geq1}[\frac 1 6]} && {\mc P_3 \mc S_*^{\geq1}[\frac 1 {24}]} \\
			&& {\mc P_1 \mc S_*^{\geq1}[\frac 1 6]} && {\mc P_2 \mc S_*^{\geq1}[\frac 1 {24}]} && \cdots \\
			&&& {\mc P_1 \mc S_*^{\geq1}[\frac 1 {24}]} && \cdots \\
			&&&& \cdots
			\arrow[from=1-1, to=2-2]
			\arrow[from=1-3, to=2-2]
			\arrow[from=1-3, to=2-4]
			\arrow[from=1-5, to=2-4]
			\arrow[from=1-5, to=2-6]
			\arrow[from=1-7, to=2-6]
			\arrow[from=2-2, to=3-3]
			\arrow[from=2-4, to=3-3]
			\arrow[from=2-4, to=3-5]
			\arrow[from=2-6, to=3-5]
			\arrow[from=2-6, to=3-7]
			\arrow[from=3-3, to=4-4]
			\arrow[from=3-5, to=4-4]
			\arrow[from=3-5, to=4-6]
			\arrow[from=4-4, to=5-5]
		\end{tikzcd}\]
		However, the inclusion of the former into the latter is ($\infty$)-initial so both definitions agree, and we decide to stick with the first one.
	\end{rmq}
There is a natural $\bb S$-tamification functor $\mc S_*^{\geq 1} \to \mc{FT}$ that takes a pointed, connected space and constructs the diagram 
\[\begin{tikzcd}[row sep=large, column sep=-10]
	{\mc P_1 X} && {\mc P_2 X[\frac 1 2]} && {\mc P_3 X[\frac 1 6]} && \cdots \\
	& {\mc P_1 X[\frac 1 2]} && {\mc P_2 X[\frac 1 6]} && {\mc P_3 X[\frac 1 {24}]}
	\arrow["{[\frac 1 2]}"', from=1-1, to=2-2]
	\arrow["{\mc P_1}", from=1-3, to=2-2]
	\arrow["{[\frac 1 3]}"', from=1-3, to=2-4]
	\arrow["{\mc P_2}", from=1-5, to=2-4]
	\arrow["{[\frac 1 4]}"', from=1-5, to=2-6]
	\arrow["{\mc P_3}", from=1-7, to=2-6]
\end{tikzcd}\]
that is, $X_n = \mc P_n X[\frac 1 {n!}]$, with the equivalences $X_n[\frac 1 {n+1}] \simeq \mc P_n X[\frac 1 {(n+1)!}] \simeq \mc P_n X_{n+1}$ induced by $\mc P_n \mc P_{n+1} \simeq \mc P_n$.
We will write this functor $(-)^{\bb S\text{-tame}}$.
\begin{rmq}
	Note that the $\bb S$-tamification functor is \emph{not} a localization. It is actually not even essentially surjective, for example 
	\[\begin{tikzcd}[row sep=large, column sep=10]
		{\Sigma H\bb F_2} && {0} && {0} && \cdots \\
		& {0} && {0} && {0}
		\arrow[from=1-1, to=2-2]
		\arrow[from=1-3, to=2-2]
		\arrow[ from=1-3, to=2-4]
		\arrow[ from=1-5, to=2-4]
		\arrow[from=1-5, to=2-6]
		\arrow[from=1-7, to=2-6]
	\end{tikzcd}\]
	is an $\bb S$-tame space which is not in the image of the tamification functor.
\end{rmq}
We can now prove \Cref{Stamemodel}. This is simply a consequence of \Cref{modelgoodwillie}. The only remaining point to check is that modulo the identification of \Cref{modelgoodwillie}, the functors $\mc P_{n-1}\colon \mc P_n \mc S_*^{\geq 1} \to \mc P_{n-1}\mc S_*^{\geq 1}$ amount to the natural map from $n$-nilpotent Lie algebras to $(n-1)$-nilpotent Lie algebras, namely the left adjoint
$$\mrm{Nil}_{n-1}\colon \mrm{Alg}(\tau_n \bb L)^{\geq 1} \to \mrm{Alg}(\tau_{n-1}\bb L)^{\geq 1}$$
to the restriction along the truncation $\tau_n \bb L \to \tau_{n-1}\bb L$, and that the localizations $[\frac 1 {k!}]$ correspond on both sides, coherently in $n$. The second point is easy. For the first point, after inverting $n!$ the Tate constructions of $\Sigma_k$ vanish for $k \leq n$, so the spaces of $k$-stages for $\mrm{Sp}^\otimes$ are contractible for $k \leq n$ \cite[Cor.~1.25]{heuts2021}. Since $\mc P_{n-1}$ sends an $n$-stage to an $(n-1)$-stage, the equivalences $\mc P_k \mc S_* \simeq \mrm{coAlg}^{\mrm{ind}}(\tau_k \mrm{Sp}^{\otimes})$ of \cite[Cor.~1.26]{heuts2021} are compatible, up to a contractible space of choices, with $\mc P_{n-1}$ on the left and with the functor $\mrm{coAlg}^{\mrm{ind}}(\tau_n \mrm{Sp}^{\otimes}) \to \mrm{coAlg}^{\mrm{ind}}(\tau_{n-1} \mrm{Sp}^{\otimes})$ forgetting the $n$-fold comultiplication on the right (see the proof of \cite[Lemma~4.31]{heuts2021}, where the Tate terms now vanish). On the algebraic side, by \cite[Prop.~6.5, 6.6 and 6.9]{heuts2021}, the functor $\mbf{bar}_{\tau_n \bb L}$ identifies the Goodwillie tower of $\mrm{Alg}(\tau_n \bb L)$, whose stages are given by the functors $\mrm{Nil}_k$, with the tower of $\infty$-categories $\mrm{coAlg}^{\mrm{ind}}_{\mrm{dp}}(B(\tau_k \bb L)) \simeq \mrm{coAlg}^{\mrm{ind}}(\tau_k \mrm{Sp}^{\otimes})$, the last equivalence holding since the relevant Tate constructions vanish. Hence the equivalences of \Cref{modelgoodwillie} for $n$ and $n-1$ intertwine $\mc P_{n-1}$ and $\mrm{Nil}_{n-1}$. {The levelwise equivalences thus assemble into an equivalence between the diagram of \Cref{formaltamespace} and that of \Cref{def:spectrallazard}, and passing to the limit gives $\mc{FT} \simeq \mrm{LazAlg}(\bb S)$.}

\section{Recovering mild and tame homotopy theory}\label{section:mildandtame}
	By suitably truncating, Lie algebra models can be defined over the integers.
	\begin{defn}\label[defn]{def0tame}
	Fix an integer $r \geq 1$ and a nondecreasing function $s: \bb N \to \bb N \cup \infty$. By convention we say $t_{\infty}$ is the identity functor. 
	We say that a map $f: X \to Y$ between $(r-1)$-connected homotopy types is an $(r,s)$-\emph{tame} equivalence if {$t_{s(n)} \mc P_{n-1} f$ is an equivalence after inverting $(n-1)!$, for every integer $n \geq 2$, where $\mc P_{n-1}$ is the $(n-1)$-th stage} of the Goodwillie tower of the identity.  We simply call $1$-tame equivalence a $(1,s)$-tame equivalence with $s$ the constant function with value  $\infty$. 
\end{defn}
	\begin{defn}
		Fix a nondecreasing function $s: \bb N \to \bb N$ and as always a parameter $r \geq 1$. 
		The category of $(r,s)$-tame spaces is the following limit of $\infty$-categories{, where at the stage $\mc P_{n-1}$, in which $(n-1)!$ is inverted, we keep the homotopy groups in degrees $[r, s(n)]$}:
		\[\begin{tikzcd}[row sep=large, column sep=-10]
			{\mc P_1 \mc S_*[\frac 1 1]^{[r, s(2)]}} && {\mc P_2 \mc S_*[\frac 1 2]^{[r,s(3)]}} && {\mc P_3 \mc S_*[\frac 1 6]^{[r,s(4)]}} && \cdots \\
			& {\mc P_1 \mc S_*[\frac 1 2]^{[r,s(2)
					]}} && {\mc P_2 \mc S_*[\frac 1 6]^{[r,s(3)]}} && {\mc P_3 \mc S_*[\frac 1 {24}]^{[r,s(4)]}}
			\arrow[from=1-1, to=2-2]
			\arrow[from=1-3, to=2-2]
			\arrow[from=1-3, to=2-4]
			\arrow[from=1-5, to=2-4]
			\arrow[from=1-5, to=2-6]
			\arrow[from=1-7, to=2-6]
		\end{tikzcd}\]

	\end{defn}
	{Let $X$ be $(r-1)$-connected. By \cite[Thm.~1.13]{goodwillie03}, the map $X \to \mc P_{n-1} X$ induces isomorphisms on $\pi_i$ for $i \leq n(r-1)$ (see also \cite[\S 6.3]{heuts2021}). We obtain the following description.} 
	\begin{lemme}\label[lemme]{lemme:lazardtronque}
		If {$s(n) \leq n(r-1)$ for all $n \geq 2$}, the category of $(r,s)$-tame spaces is equivalent to the limit of the following diagram.
				\[\begin{tikzcd}[row sep=large, column sep=-10]
			{\mc S_*[\frac 1 1]^{[r, s(2)]}} && { \mc S_*[\frac 1 2]^{[r,s(3)]}} && { \mc S_*[\frac 1 6]^{[r,s(4)]}} && \cdots \\
			& { \mc S_*[\frac 1 2]^{[r,s(2)
					]}} && { \mc S_*[\frac 1 6]^{[r,s(3)]}} && {\mc S_*[\frac 1 {24}]^{[r,s(4)]}}
			\arrow[from=1-1, to=2-2]
			\arrow[from=1-3, to=2-2]
			\arrow[from=1-3, to=2-4]
			\arrow[from=1-5, to=2-4]
			\arrow[from=1-5, to=2-6]
			\arrow[from=1-7, to=2-6]
		\end{tikzcd}\]
	Moreover, this category is equivalent to the category of $(r,s)$-tame spectral Lazard algebras defined by the limit
\[\begin{tikzcd}[row sep=large, column sep=-30]
	{\mrm{Alg}(\tau_{1} \bb L)[\frac 1 1]^{[r, s(2)]}} && {\mrm{Alg}(\tau_{2} \bb L)[\frac 1 2]^{[r,s(3)]}} && {\mrm{Alg}(\tau_{3} \bb L)[\frac 1 6]^{[r,s(4)]}} && \cdots \\
	& {\mrm{Alg}(\tau_{1} \bb L)[\frac 1 2]^{[r,s(2)]}} && {\mrm{Alg}(\tau_{2} \bb L)[\frac 1 6]^{[r,s(3)]}} && {\mrm{Alg}(\tau_{3} \bb L)[\frac 1 {24}]^{[r,s(4)]}}
	\arrow[from=1-1, to=2-2]
	\arrow[from=1-3, to=2-2]
	\arrow[from=1-3, to=2-4]
	\arrow[from=1-5, to=2-4]
	\arrow[from=1-5, to=2-6]
	\arrow[from=1-7, to=2-6]
\end{tikzcd}\]
Here $\mrm{Alg}(\tau_k \bb L)[\frac1m]^{[a,b]}$ denotes the $\infty$-category of $\tau_k\bb L$-algebras whose underlying spectrum is $m$-local with homotopy groups concentrated in degrees $[a,b]$.
	\end{lemme}
	\begin{lemme}\label[lemme]{lemme:simpliflazardtronque}
		Suppose $r \geq2$ and {$s(n) \leq n(r-1)$ for all $n \geq 2$}.
		\begin{enumerate}
			\item The category of $(r,s)$-tame spaces is equivalent to the full subcategory of $(r-1)$-connected spaces $X$ such that $p$ is invertible in {$\pi_k X$ if $k > s(p)$}.
			\item There is an equivalence between $(r,s)$-tame spectral Lazard algebras and {$(r-1)$-connected} spectral Lie algebras $L$ such that $p$ is invertible in {$\pi_k L$ if $k > s(p)$}.
			\item Moreover, if {$s(n) \leq r+2n-4$ for all $n \geq 2$}, then $(r,s)$-tame spectral Lazard algebras can simply be described by $(r-1)$-connected algebras over the integral Lie operad $\bb L(\bb Z)$ with the same divisibility conditions on the {$\pi_k$}'s.
		\end{enumerate}

	\end{lemme}
	\begin{proof}
		The point $1)$ is straightforward. For the second point, one observes that {$$\mrm {Alg}(\tau_{n-1} \bb L)^{[r, s(n)]} \simeq \mrm{Alg}(\bb L)^{[r, s(n)]}$$} by \cite[proof of Thm.~0.4]{heuts2021}. This is simply for the reason that {the $n$-th term in the lower central series of an $(r-1)$-connected algebra is $n(r-1)$-connected}. Finally the third point follows as in \Cref{modelgoodwillieZ}, see \cite[Lemma 3.7.13]{wu22}. This is essentially because $p$-locally, $t_{2p-4} \bb S \simeq H\bb Z_{(p)}${, and at the stage $\mc P_{n-1}$ the first prime which is not inverted is at least $n$}.
	\end{proof}
Now \Cref{2tame} and \Cref{mildmodel} are easy consequences of \Cref{lemme:lazardtronque} together with \Cref{lemme:simpliflazardtronque}.

\section{What $\bb S$-tamification remembers about the fundamental group}\label{section:pi1}
	Note that due to convergence issues, the $\bb S$-tame homotopy type of a non-simply connected space does not necessarily contain the information of the rational homotopy type, contrary to the $(r-1)$-connected theory with $r \geq 2$. It allows one to recover, though, the information of the Lazard completion, and thus the Malcev completion of the fundamental group. The purpose of this section is to prove this claim.
	\begin{prop}
		Let $X$ be a pointed connected space. Its $\bb S$-tamification allows one to recover the Lazard completion of $\pi_1 X$. 
	\end{prop}
	This is a consequence of the following folklore proposition: 

	\begin{prop}\label[prop]{pi1goodwillie}
		Let $X$ be a pointed connected space. The map $p_n\colon X\to\mc P_nX$ induces an isomorphism
		\[\pi_1X/\Gamma^{n+1}\pi_1X\xrightarrow{\ \cong\ }\pi_1\mc P_nX,\]
		where $\Gamma^1G=G$ and $\Gamma^{k+1}G=[G,\Gamma^kG]$ is the lower central series.
	\end{prop}
	\begin{proof}
		For a group $H$ write $N_nH=H/\Gamma^{n+1}H$. All colimits of spaces are pointed homotopy colimits.
		
		Firstly, we claim that on pointed connected spaces, $\mc P_n$ takes connected values and commutes with
			filtered colimits and geometric realizations.
		For filtered colimits, $\mc P_n$ is built from finite homotopy limits and sequential colimits
		\cite[\S1]{goodwillie03}. For the rest we argue by induction on $n$. For $n=1$,
		$\mc P_1=\Omega^\infty\Sigma^\infty$, and $\Omega^\infty$ commutes with sifted colimits of
		connective spectra \cite[1.4.3.9]{Lurie2015}. For $n\geq2$, by \cite[Lemma 2.2]{goodwillie03} there is
		a natural fiber sequence
		\[\mc P_nX\to\mc P_{n-1}X\to R_nX,\qquad R_nX=\Omega^\infty E_n(X),\quad
		E_n(X)=\Sigma\,\partial_n(\mathrm{id})\wedge_{h\Sigma_n}(\Sigma^\infty X)^{\wedge n}.\]
		As $\partial_n(\mathrm{id})$ is non-equivariantly a wedge of copies of $S^{1-n}$
		\cite{johnson95,aronemahowald}, the spectrum $E_n(X)$ is $1$-connected for $X$ connected. Hence
		$R_nX$ is simply connected and $\mc P_nX$ is connected by the long exact sequence. The functor
		$E_n$ commutes with sifted colimits, hence so does
		$R_n$ on connected spaces. Let $X_\bullet$ be a simplicial pointed connected space. Since
		$\mc P_{n-1}X_\bullet$ and $R_nX_\bullet$ are levelwise connected, realization commutes with taking the homotopy fiber $\mc P_nX_\bullet \simeq \mrm{fib}( \mc P_{n-1}X_\bullet \to {R_nX_\bullet})$
		\cite[Theorem B.4]{bousfieldfriedlander}, so the claim is true for $\mc P_n$.
		Secondly, we reduce the statement to the case where $X$ is a wedge of circles.
		The functor $\pi_1$ from pointed connected spaces to groups is left adjoint to $B$, so it turns
		homotopy colimits into colimits. The functor $N_n$ commutes with filtered colimits and
		coequalizers, since $n$-nilpotent groups are closed under these in groups. By the first claim, the
		functors $\pi_1\mc P_n$, $N_n\pi_1$ and $N_n\pi_1\mc P_n$ turn filtered colimits and geometric
		realizations of pointed connected spaces into colimits of groups. Therefore the class $\mathcal C$
		of pointed connected spaces $X$ for which both natural maps
		\[q_X\colon\pi_1\mc P_nX\to N_n\pi_1\mc P_nX,\qquad
		N_n\pi_1(p_n)\colon N_n\pi_1X\to N_n\pi_1\mc P_nX\]
		are isomorphisms is closed under filtered colimits and geometric realizations.
		
		It remains to prove the statement for $X=\bigvee_kS^1$ with $k$ finite. By \cite[Theorem 6.8]{biedermanndwyer}, for the free group $F_k=\pi_1X$,
		$\pi_1(p_n)\colon F_k\to\pi_1\mc P_nX$ is surjective with kernel $\Gamma^{n+1}F_k$. Thus
		$\pi_1\mc P_nX$ is $n$-nilpotent, $q_X$ is an isomorphism, and so is $N_n\pi_1(p_n)$; that is,
		$X\in\mathcal C$. Arbitrary wedges of circles are filtered colimits of finite ones, hence lie in
		$\mathcal C$. Finally, a pointed connected space is modeled by a reduced simplicial set $X$. Let $GX$ be
Kan's loop group of $X$ (it is a simplicial group with each $G_mX$ free, modeling $\Omega X$) and $\overline W$ the
classifying space functor for simplicial groups; the unit $X\to\overline WGX$ is a weak
equivalence \cite[Ch.~V]{goerssjardine}.   Hence $X\simeq|[m]\mapsto
B(G_mX)|$ is a geometric realization of wedges of circles and $X\in\mathcal C$. The isomorphism of the statement is $q_X^{-1}\circ N_n\pi_1(p_n)$.
	\end{proof}
Now let $(X_n)$ be an $\bb S$-tame homotopy type coming from a pointed space $X$, represented by a zigzag diagram with the coherences {$\mc P_{n-1} X_n \simeq X_{n-1} [\frac 1 {n!}]$}. A simple consequence of the previous proposition is that applying the $\pi_1$ functor, we get a diagram in groups
\[\begin{tikzcd}[row sep=large, column sep=10]
	{G/\Gamma^2} && {(G/\Gamma^3)[\tfrac12]} && {(G/\Gamma^4)[\tfrac16]} && \cdots \\
	& {(G/\Gamma^2)[\tfrac12]} && {(G/\Gamma^3)[\tfrac16]} && {(G/\Gamma^4)[\tfrac1{24}]}
	\arrow[from=1-1, to=2-2]
	\arrow[two heads, from=1-3, to=2-2]
	\arrow[from=1-3, to=2-4]
	\arrow[two heads, from=1-5, to=2-4]
	\arrow[from=1-5, to=2-6]
	\arrow[two heads, from=1-7, to=2-6]
\end{tikzcd}\]
where $G = \pi_1 X$. {By \Cref{equivalencelazarddefs}, transported to groups through Lazard's correspondence \cite[Chap.~II, \S 4]{lazard54}, this is exactly the data of the Lazard completion of $\pi_1 X$.}

\bibliographystyle{alpha}
\bibliography{tame_goodwillie}
\end{document}